\documentclass[a4paper,12pt]{article}
\usepackage[T2A]{fontenc}                      
\usepackage[cp1251]{inputenc}           
\usepackage[all]{xy}
\usepackage{amssymb}
\usepackage{cite}
\usepackage{cmap}
\usepackage{latexsym}
\usepackage{enumerate}
\usepackage{amsmath, amsthm, amscd, amsfonts, amssymb, graphicx, color}
\usepackage[left=2.5cm,right=2.5cm,top=2cm,bottom=2cm]{geometry}
\usepackage{indentfirst}
\usepackage{array}
\usepackage{bm}
\usepackage{float}

\usepackage{wrapfig}

\usepackage{authblk}
\usepackage{enumitem}

\newtheorem{proposition}{Proposition}

\newtheorem{definition}{Definition}

\newtheorem{theorem}{Theorem}     
\newtheorem{lemma}{Lemma}         
\newtheorem{corollary}{Corollary} 
\newtheorem{conjecture}{Conjecture}

\usepackage[colorlinks=true, linkcolor=blue, citecolor=blue, urlcolor=blue]{hyperref}

\title{Every tournament $T$ of order at least 3\\  admits infinitely many $T$-irregular oriented graphs}

\author{Tatiana~Dovzhenok}

\affil{\small Research Laboratory ``Mathematics of Hybrid Intelligence Systems'', \\
	Francisk Skorina Gomel State University, Gomel, 246028, Belarus}

\date{} 

\begin{document}	
\maketitle

\begin{abstract}
Let $F$ be a fixed graph. A graph $G$ is $F$-irregular if any two distinct vertices of $G$ belong to a different number of subgraphs of $G$ isomorphic to $F$. In this paper, we settle the strong conjecture about $F$-irregular oriented graphs (Dovzhenok, Lukashenko, and Filiuta, 2026) for the case of tournaments. Specifically, we prove that for every tournament $T$ of order $n \ge 3$, there exists an infinite family of \mbox{$T$-irregular} oriented graphs. Our proof is constructive. For $n \ge 4$, it is based on the vertex cloning method, and the underlying graphs of the constructed $T$-irregular oriented graphs are inherently $K_n$-irregular.
\smallskip
	
\end{abstract}

\section{Introduction}
For a fixed graph $F$, the \emph{$F$-degree} of a vertex $v$ in a graph $G$ is the number of subgraphs of $G$ isomorphic to $F$ containing $v$~\cite{1, 2}. A graph $G$ is \emph{$F$-irregular} if all its vertices have pairwise distinct $F$-degrees. Such a generalization extends the traditional notion of vertex degrees to higher-order structural motifs, resolving the classical limitation that no simple finite graph can have vertices of pairwise distinct degrees. This line of research was motivated by the following conjecture.

\begin{conjecture}[Chartrand, Holbert, Oellermann, and Swart, 1987~\cite{1}] \label{con1}
	There exists a~non-trivial $F$-irregular graph for every connected graph $F$ of order at least $3$.
\end{conjecture}

The accumulation of positive results for various graph classes suggested that \mbox{$F$-irregular} graphs are not merely isolated exceptions, but exist in abundance.  This shift from existence to cardinality naturally extended the original inquiry.

\begin{conjecture} 
	[Strong conjecture about $F$-irregular graphs, Dovzhenok, Filuta, and Chuhai, 2024~\cite{3}] 
	\label{con2}
	For every connected graph $F$ of order $|F| \ge 3$, there exist infinitely many $F$-irregular graphs.
\end{conjecture}

The authors confirmed this conjecture for biconnected graphs with minimum degree~$2$~\cite{3}, while the first author settled the case for graphs $F$ of diameter $2$~\cite{4}. A series of recent preprints exhausts all cases. Specifically, Schreib~\cite{5} proposed an algorithm for constructing $F$-irregular graphs for non-complete connected patterns $F$. Independently, Dovzhenok and Filuta~\cite{6} provided a complete proof of Conjecture~\ref{con2}. Finally, a solution for arbitrary disconnected $F$ was obtained by Dovzhenok in~\cite{7}.

In parallel, the concept of $F$-irregularity was carried over to oriented graphs in~\cite{8}, where an oriented analogue of Conjecture~\ref{con2} was formally introduced and explicitly verified for all directed cycles $\overrightarrow{C_n}$ with $n \ge 3$.

\begin{conjecture} 
	[Strong conjecture about $F$-irregular oriented graphs, Dovzhenok, Luka\-shen\-ko, and Filiuta, 2026~\cite{8}] \label{con3}
	For every weakly connected oriented graph $F$ of order $|F| \ge 3$, there are infinitely many $F$-irregular oriented graphs.
\end{conjecture}

In this paper, we settle Conjecture~\ref{con3} for the class of tournaments. Recall that a \emph{tournament} $T$ of order $n$ is an oriented graph obtained by assigning a direction to each edge of a complete graph $K_n$. Our main result states the following:

\begin{theorem}\label{t1}
	Every tournament $T$ of order $n \ge 3$ admits an infinite family of $T$-irregular oriented graphs.
\end{theorem}

The proof of Theorem~\ref{t1} is split into two parts. Section~\ref{sec3} investigates the base case of tournaments of order $3$. Section~\ref{sec4} presents a general approach to constructing $T$-irregular oriented graphs for $n \ge 4$ via the vertex cloning method. These constructions enable the exact computation and comparison of all vertex $T$-degrees.
 
Furthermore, for every tournament $T$ of order $n \ge 4$, the underlying graphs of the \mbox{$T$-irregular} oriented graphs constructed in Section~\ref{sec4} are inherently $K_n$-irregular (Corollary~\ref{cor.5}). This provides an alternative proof of Conjecture~\ref{con2} for complete graphs $K_n$ of order $n \ge 4$ to that given in~\cite{6}.

To conclude this section, we note that the problem of $F$-irregularity for general digraphs trivially reduces to the undirected case. Indeed, if $U(F)$ denotes the underlying graph of a weakly connected digraph $F$, any $U(F)$-irregular graph can be transformed into an $F$-irregular digraph by replacing each edge with a pair of symmetric arcs. Consequently, the oriented setting remains the primary non-trivial domain for this problem.

\section{Preliminaries} \label{sec2}
Unless stated otherwise, all numerical variables throughout this paper are assumed to be positive integers.

The primary objects of our study are finite oriented graphs. An oriented graph $G$ is obtained by assigning a direction to each edge of a simple undirected graph. Let $V(G)$ and $A(G)$ denote the vertex set and the arc set of $G$, respectively. An arc directed from a vertex $u$ to a vertex $v$ is represented by the ordered pair $(u, v)$. The cardinality $|V(G)|$, written more concisely as $|G|$, is the \emph{order} of $G$, and $G$ is \emph{non-trivial} if $|G| \ge 2$. An \emph{independent set} in $G$ is a subset $I$ of $V(G)$ such that no two vertices in $I$ are adjacent. 

For any graphs $F$ and $G$ (both either undirected or oriented), an \emph{$F$-subgraph} of $G$ is a subgraph of $G$ isomorphic to $F$. For an oriented graph $G$, a subgraph $H$ of $G$ is \emph{induced} if, for any pair of vertices $u, v \in V(H)$, the arc $(u, v)$ belongs to $A(H)$ if and only if $(u, v) \in A(G)$.
The \emph{underlying graph} of an oriented graph $G$, denoted by $U(G)$, is the undirected graph obtained by replacing each directed arc $(u, v) \in A(G)$ with an undirected edge $uv$. An oriented graph $G$ is \emph{weakly connected} if $U(G)$ is connected. 
\vspace{-1mm}
\section{Tournaments of order $3$} \label{sec3}
We first analyze the base case of tournaments $T$ of order $3$. Up to isomorphism, there are only two such tournaments: the \emph{directed 3-cycle} $\overrightarrow{C_3}$ with arcs $(v_1, v_2)$, $(v_2, v_3)$, and $(v_3, v_1)$, and the \emph{transitive tournament} $TT_3$ with arcs $(v_1, v_2)$, $(v_2, v_3)$, and $(v_1, v_3)$, both defined on the vertex set $\{v_1, v_2, v_3\}$.

\begin{proposition}\label{p1}
	For every tournament~$T$ of order $3$, there exist infinitely many $T$-irregular oriented graphs.
\end{proposition}

\begin{proof}
	Let $T$ be a tournament of order $3$. We consider two cases depending on $T$.
	
	\medskip
	\noindent\textbf{Case 1.} $T \cong \overrightarrow{C_3}$.
	The result follows immediately from \cite{8}.
	
	\medskip
	\noindent\textbf{Case 2.} $T \cong TT_3$.
	By \cite{3}, there exists an infinite family $\mathcal{G}$ of finite, simple, undirected \mbox{$K_3$-irregular} graphs. For any $H \in \mathcal{G}$, we construct an oriented graph $G$ by assigning an acyclic orientation to $H$, obtained by fixing a total ordering on $V(H)$ and directing all edges from lower to higher indices. Since $G$ is acyclic, it contains no directed \mbox{$3$-cycles}. Consequently, a triplet of vertices forms a $TT_3$-subgraph in $G$ if and only if it forms a $K_3$-subgraph in $H$. This bijective correspondence ensures that for every $v \in V(G)$, the number of $TT_3$-subgraphs containing $v$ in $G$ equals the number of $K_3$-subgraphs containing $v$ in $H$. Since $H$ is $K_3$-irregular, it follows that $G$ is $TT_3$-irregular.
\end{proof}

\section{The general case: tournaments of order $n \ge 4$} \label{sec4}
In this section, for an arbitrary fixed tournament $T$ of order $n \ge 4$, we present a constructive method to obtain an infinite family of $T$-irregular oriented graphs. The proposed graphs $D(T, \bar{p})$ rely on auxiliary tuples of prime numbers $\bar{p} = (p_1, p_2, \dots, p_{n-1})$. We begin by introducing a specific number-theoretic condition on these primes, which  enables us to prove that the $T$-degrees of the vertices in $D(T, \bar{p})$ are pairwise distinct.
\subsection{Main Condition}

\begin{definition}
	Let $n \ge 4$. An ordered tuple of prime numbers $\bar{p} =(p_1, p_2, \dots, p_{n-1})$ satisfies the \textbf{main condition} if 
	\begin{equation} \label{eq:main}
		p_1 > \max(6, n), \qquad p_{i+1} > 3p_i \quad \text{for all } 1 \le i \le n-2, \qquad \text{and} \qquad p_{n-1} < 6p_{n-2}.
	\end{equation}
\end{definition}

\begin{lemma} \label{l1}
	For any $n \ge 4$, there exist infinitely many ordered tuples of prime numbers $\bar{p} =(p_1, p_2, \dots, p_{n-1})$ that satisfy the main condition and have distinct sums.
\end{lemma}
\begin{proof}
	By Euclid's theorem, there exist infinitely many prime numbers. Let us choose an initial prime $p_1 > \max(6, n)$. Successively, we select primes $p_2, p_3, \dots, p_{n-2}$ such that $p_{i+1} > 3p_i$ for all $1 \le i \le n-3$. Finally, by Bertrand's postulate, there exists a prime $p_{n-1}$ such that $3p_{n-2} < p_{n-1} < 6p_{n-2}$. Thus, we obtain a tuple satisfying the main condition. Since $p_1$ can be chosen arbitrarily large, the resulting sums $S = \sum_{j=1}^{n-1} p_j$ can also be made arbitrarily large, ensuring infinitely many such tuples with distinct sums.
\end{proof}

\begin{lemma} \label{l2}
	Let $n \ge 4$ and let $\bar{p} =(p_1, p_2, \dots, p_{n-1})$ be a tuple of prime numbers conforming to \eqref{eq:main}. Then the product $P = \prod_{j=1}^{n-1} p_j$ and the sum $S = \sum_{j=1}^{n-1} p_j$ satisfy
	\[
	P > p_i p_{n-1} \quad \text{for all } 1 \le i \le n-1 \qquad \text{and} \qquad S > 3p_i + n \quad \text{for all } 1 \le i \le n-2.
	\]
\end{lemma}

\begin{proof}
	By the main condition \eqref{eq:main}, we have $p_1 > 6$ and $6p_{n-2} > p_{n-1}$. Since $p_{i+1} > 3p_i$, the sequence of primes is strictly increasing, which implies $p_{n-1} > 3p_{n-2} \ge 3p_i$ and $p_i \ge p_1 > n$ for all $1 \le i \le n-2$. Given that $n \ge 4$, the desired inequalities follow directly from the estimates:
	\begin{align*}
		P &\ge p_1 p_{n-2} p_{n-1} > 6 p_{n-2} p_{n-1} > p_{n-1}^2 \ge p_i p_{n-1} \quad \text{for all } 1 \le i \le n-1, \\[1mm]
		S &> p_i + p_{n-1} > n + 3p_i \quad \text{for all } 1 \le i \le n-2. \qedhere
	\end{align*}
\end{proof}

\subsection{The Graph $D(T, \bar{p})$}
\subsubsection{Vertex cloning method}
To obtain our family of $T$-irregular oriented graphs, we employ the vertex cloning method. Under this approach, each vertex of $D(T, \bar{p})$ is a \emph{clone} of a specific vertex in the tournament $T$, called its \emph{prototype}. Clones with identical prototypes are non-adjacent in $D(T, \bar{p})$.
If two clones correspond to distinct prototypes in $T$ and are adjacent in $D(T, \bar{p})$, their connecting arc is oriented consistently with the arc between their prototypes in $T$. Formally, let $u_0$ and $v_0$ be distinct vertices of $T$, and let $u$ and $v$ be clones of $u_0$ and $v_0$, respectively. We denote the arc between the clones $u$ and $v$ in $D(T, \bar{p})$ by $(u, v)^*$, defined as follows:
\[
(u, v)^* = 
\begin{cases} 
	(u, v), & \text{if } (u_0, v_0) \in A(T); \\ 
	(v, u), & \text{if } (v_0, u_0) \in A(T). 
\end{cases}
\]
By analogy with standard graph notation, the clones $u$ and $v$ are adjacent in $D(T, \bar{p})$ if and only if $(u, v)^* \in A(D(T, \bar{p}))$.

\subsubsection{Structure of the Graph $D(T, \bar{p})$}
Without loss of generality, we assume that $V(T) = \{1, 2, \dots, n\}$. 
\smallskip

Consider an arbitrary ordered tuple of prime numbers $\bar{p} = (p_1, p_2, \dots, p_{n-1})$ satisfying the main condition \eqref{eq:main}. Let $S = \sum_{j=1}^{n-1} p_j$ and $P = \prod_{j=1}^{n-1} p_j$.
\smallskip

We define the graph $D(T, \bar{p})$ of order $2S-n+2$ via its vertex set  
\[
V(D(T, \bar{p})) = \left( \bigcup_{i=1}^{n-1} A_i \right) \cup \left( \bigcup_{i=1}^{n-1} B_i \right) \cup \{c\}, 
\]
where the constituent vertex subsets are given by
\begin{align*}
A_i &= \bigl\{ a_{i,j} \bigm| 1 \le j \le p_i, \text{ and } a_{i,j} \text{ is a clone of vertex } i \bigr\},\\[1mm] 
B_i &= \bigl\{ b_{i,j} \bigm| 1 \le j \le p_i - 1, \text{ and } b_{i,j} \text{ is a clone of vertex } n \bigr\},\\[1mm]
c &\quad \text{is a clone of vertex } n;
\end{align*}
and its arc set 
\[
A\bigl(D(T, \bar{p})\bigr) = \left( \bigcup_{i=1}^{n-1} E_i \right) \cup \left( \bigcup_{i=1}^{n-2} \bigcup_{k=i+1}^{n-1} E_{i,k} \right) \cup X, 
\]
where the constituent arc subsets are specified by
\begin{align*}
	E_i &= \bigl\{ (a_{i,j}, b_{i,m})^* \bigm| 1 \le m<j \le p_i\bigr\}, \\[1mm]
	E_{i,k} &= \bigl\{ (a_{i,j}, a_{k,l})^* \bigm| 1 \le j \le p_i, \ 1 \le l \le p_k \bigr\} \\[1mm] 
	&\quad \cup \bigl\{ (a_{i,j}, b_{k,m})^* \bigm| 1 \le j \le p_i, \ 1 \le m \le p_k - 1 \bigr\} \\[1mm]
	&\quad \cup \bigl\{ (b_{i,j}, a_{k,m})^* \bigm| 1 \le j \le p_i-1, \ 1 \le m \le p_k \bigr\}, \\[1mm]
	X &= \bigl\{ (a_{i,1}, c)^* \bigm| 1 \le i \le n - 2 \bigr\} \cup \bigl\{ (a_{n-1,j}, c)^* \bigm| 1 \le j \le p_{n-1} \bigr\}.
\end{align*}

\begin{figure}[ht]
	\centering
	\includegraphics[width=10cm]{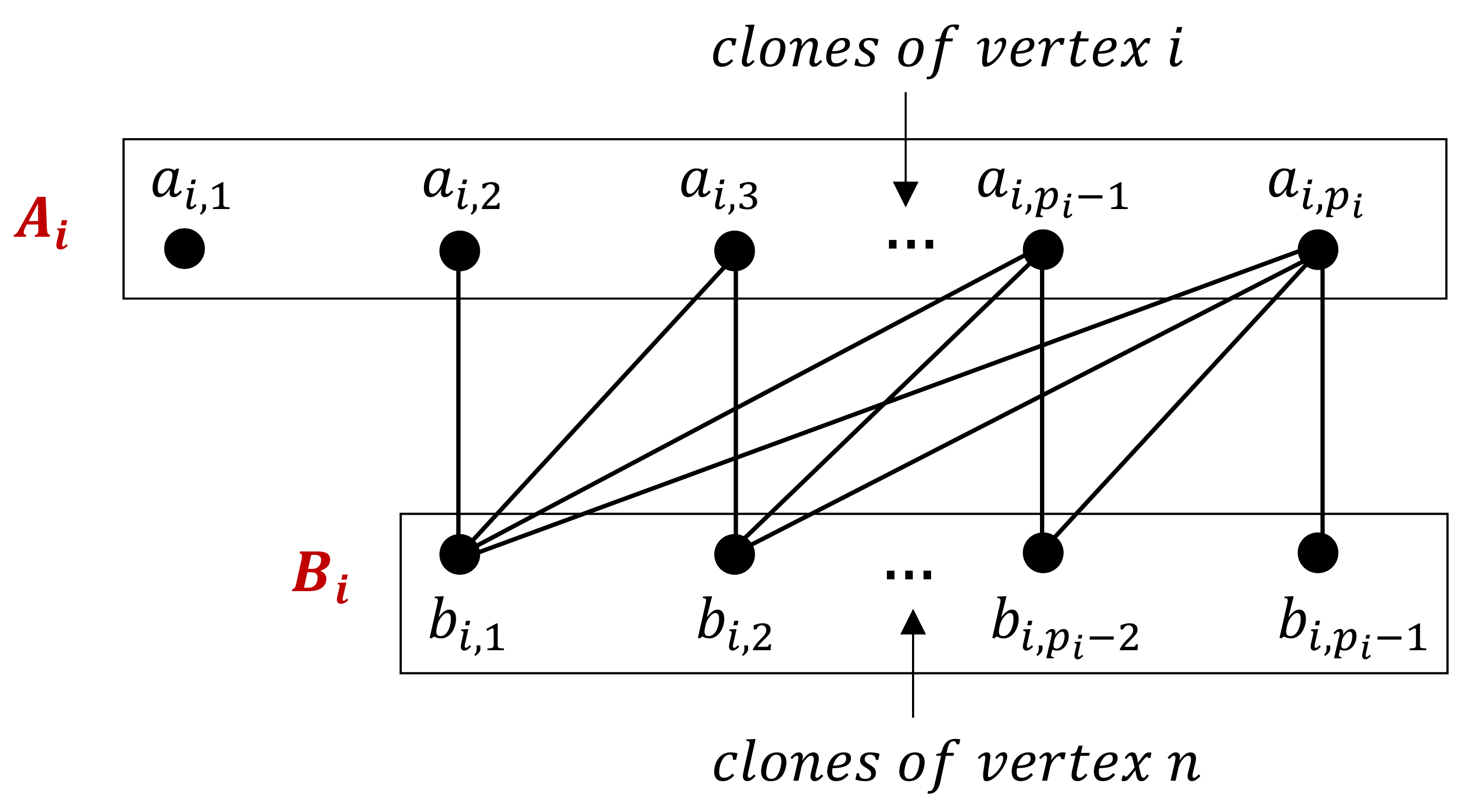}
	\caption{\small Structure of block~$i$ in~$D(T, \bar{p})$. Arcs between clones are oriented consistently with their prototypes~$i$ and $n$ in~$T$.}
	\label{fig1}
\end{figure}

For each $1 \le i \le n-1$, the vertices of $A_i$ and $B_i$ are grouped into a \emph{block} $i$ and arranged on two levels, as shown in Figure~\ref{fig1}. The upper level contains the vertices $a_{i,1}, a_{i,2}, \dots, a_{i,p_i}$ ordered from left to right, while the lower level comprises $b_{i,1}, b_{i,2}, \dots, b_{i,p_i-1}$, shifted by one position to the right. 

\begin{figure}[ht]
	\centering
	\includegraphics[width=13cm]{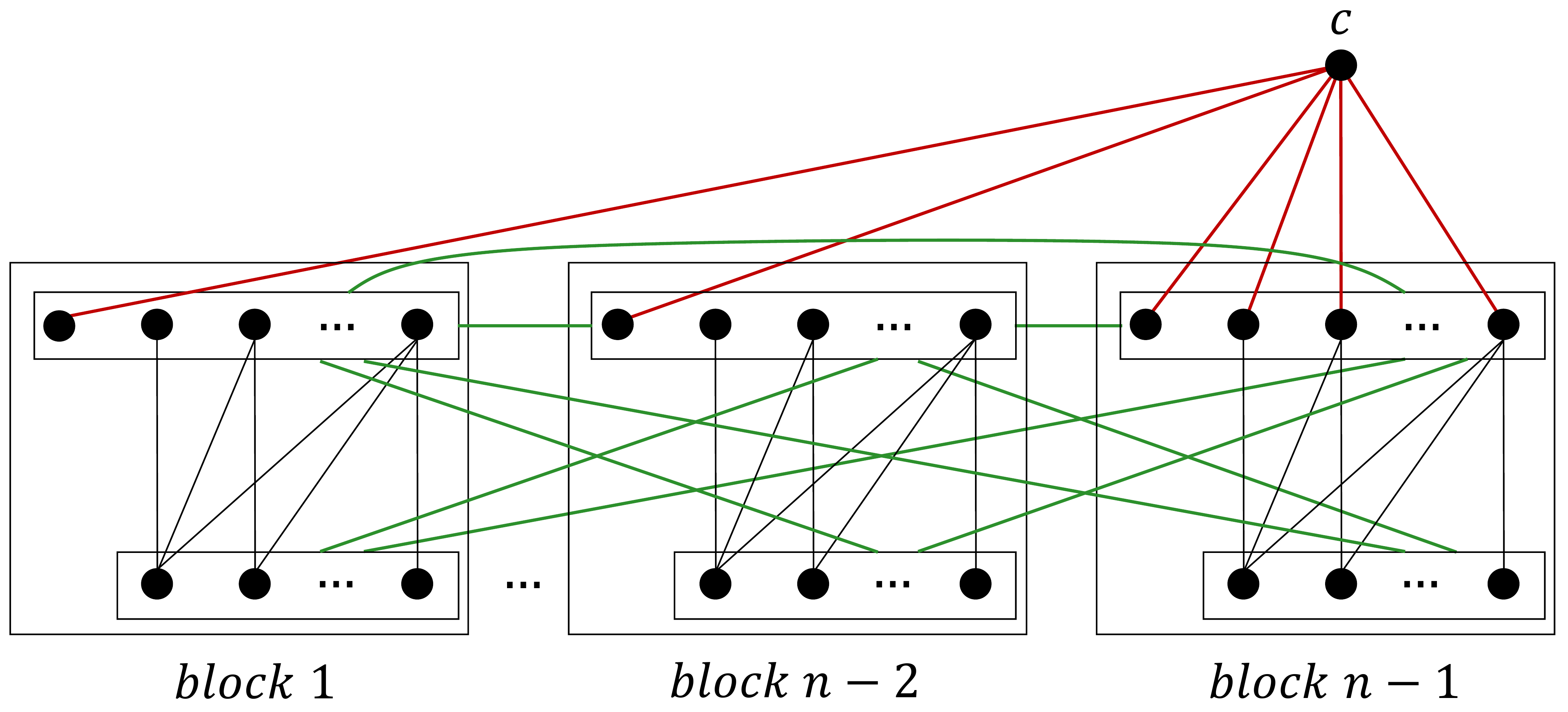}
	\caption{\small Schematic representation of the graph $D(T, \bar{p})$. Arcs between clones are oriented consistently with their prototypes in~$T$.}
	\label{fig2}
\end{figure}

Within block $i$, each lower-level vertex is joined to every upper-level vertex situated strictly above it or to its right (corresponding to the arc subset $E_i$), with no other internal arcs present.
The global structure of $D(T, \bar{p})$ is depicted in Figure~\ref{fig2}. The vertices of block $i$ are joined to those of block $s$ via the arc subsets $E_{i,s}$ (if $s > i$) or $E_{s,i}$ (if $s < i$). Specifically, every upper-level vertex of block $i$ is connected to all vertices of every other block $s \neq i$, while every lower-level vertex of block $i$ is connected to all upper-level vertices of every other block $s \neq i$. Finally, vertex $c$ is adjacent to all upper-level vertices of block $n-1$ and to the leftmost upper-level vertex $a_{i,1}$ of each remaining block $i \le n-2$ (the arc subset $X$).
\vspace{-4mm}

In this section, all $T$-degrees are considered exclusively in the graph $D(T, \bar{p})$. We denote the $T$-degrees of vertices $c$, $a_{i,j}$, and $b_{i,m}$ by $f(c)$, $f(i,j)$, and $g(i,m)$, respectively.

\subsubsection{Properties of the Graph $D(T, \bar{p})$}
\begin{lemma} \label{l3}
	A subgraph $H$ of $D(T, \bar{p})$ is a $T$-subgraph if and only if its vertices can be indexed as $V(H) = \{v_1, v_2, \dots, v_n\}$ such that the following conditions hold:
	\begin{enumerate}
		\item $H$ is an induced subgraph of $D(T, \bar{p})$;
		\item $v_i \in A_i$ for all $1 \le i \le n-1$, and $v_n \in \left( \bigcup_{l=1}^{n-1} B_l \right) \cup \{c\}$;
		\item $(v_i, v_n)^* \in A(D(T, \bar{p}))$ for all $1 \le i \le n-1$.
	\end{enumerate}
\end{lemma}

\begin{proof}
	\textbf{Necessity.} Let $H$ be a $T$-subgraph of $D(T, \bar{p})$. Since $H$ is a tournament, condition~1 is satisfied directly. To prove condition~2, note that $V(D(T, \bar{p}))$ is partitioned into $n$ independent sets: $A_1, A_2, \dots, A_{n-1}$ and $\left(\bigcup_{l=1}^{n-1} B_l \right) \cup \{c\}$. As $H$ contains no non-adjacent pairs of vertices and $|H|=n$, $H$ must have exactly one vertex in each of these sets; therefore, we can write $V(H) = \{v_1, v_2, \dots, v_n\}$ in accordance with condition~2. Finally, condition~3 follows immediately from the pairwise adjacency of all vertices in $H$.
	\smallskip
	
	\textbf{Sufficiency.} Let $H$ be a subgraph of $D(T, \bar{p})$ with $V(H) = \{v_1, v_2, \dots, v_n\}$ satisfying conditions 1--3. By the structure of $D(T, \bar{p})$, any two distinct vertices in $H$ are connected by exactly one arc. Define a mapping $\varphi: V(T) \to V(H)$ by $\varphi(i) = v_i$ for all $1 \le i \le n$. In view of condition~2, each $v_i$ is a clone of vertex $i \in V(T)$. Therefore, for any $i \neq j$, the arc $(v_i, v_j)^*$ in $H$ has the same orientation as the arc between their prototypes $i$ and $j$ in $T$. Thus, $\varphi$ establishes an isomorphism, confirming that $H$ is a $T$-subgraph of $D(T, \bar{p})$.
\end{proof}

\begin{corollary} \label{cor.1}
	The $T$-degree of a vertex $v$ in $D(T, \bar{p})$ equals the number of vertex subsets $\{v_1, v_2, \dots, v_n\}$ of $D(T, \bar{p})$ containing $v$ that satisfy conditions~2 and~3 of Lemma~\ref{l3}.
\end{corollary}

\begin{lemma} \label{l4}
	The $T$-degree of the vertex $c$ is given by
	\[
	f(c) = p_{n-1}.
	\]
\end{lemma}
\begin{proof}
	By Corollary~\ref{cor.1}, $f(c)$ equals the number of vertex subsets $\{v_1, v_2, \dots, v_{n-1}, c\}$ in $D(T, \bar{p})$ such that $v_i \in A_i$ and $(v_i, c)^* \in A(D(T, \bar{p}))$ for all $1 \le i \le n-1$. 
	\smallskip
	
	By construction, the vertex $c$ is adjacent precisely to the vertices in the set 	
	\[
	\{a_{1,1}, a_{2,1}, \dots, a_{n-2,1}\} \cup A_{n-1}.
	\] 
	This forces $v_i = a_{i,1}$ for all $1 \le i \le n-2$, while the vertex $v_{n-1} \in A_{n-1}$ can be chosen in exactly $p_{n-1}$ ways, yielding $f(c) = p_{n-1}$. 
\end{proof}

\begin{lemma} \label{l5}
	The $T$-degrees of the vertices from $B_i$ are expressed as
	\[
	g(i, m) = \frac{P}{p_i} (p_i - m) \quad \text{for all } 1 \le i \le n-1 \text{ and } 1 \le m \le p_i - 1.
	\]
\end{lemma}
\begin{proof}
	Fix $1 \le i \le n-1$ and $1 \le m \le p_i-1$. Since the vertex $b_{i,m}$ belongs to $B_i$, by Corollary~\ref{cor.1}, $g(i, m)$ equals the number of ways to choose exactly one vertex from each set $A_1, A_2, \dots, A_{n-1}$ such that $b_{i,m}$ is adjacent to all chosen vertices.
	
	First, consider the set $A_i$. In $D(T, \bar{p})$, the vertex $b_{i,m}$ is adjacent to exactly $p_i - m$ vertices in $A_i$, namely $a_{i,m+1}, a_{i,m+2}, \dots, a_{i,p_i}$, providing exactly $p_i - m$ choices. Next, for any $1 \le j \le n-1$ with $j \neq i$, the vertex $b_{i,m}$ is adjacent to every vertex in $A_j$, yielding exactly $p_j$ choices for each such set. By the product rule, the $T$-degree of $b_{i,m}$ is given by
	\[
	g(i, m) = (p_i - m) \prod_{\substack{j=1 \\ j \neq i}}^{n-1} p_j = \frac{P}{p_i} (p_i - m). \qedhere
	\] 
\end{proof}

\begin{corollary} \label{cor.2}
	$g(i, 1) > g(i, 2) > \dots > g(i, p_i-1) \quad \text{for all } 1 \le i \le n-1.$
\end{corollary}

\begin{lemma} \label{l6}
	The $T$-degrees of any vertices from distinct sets $B_i$ and $B_j$ are unequal:
	\begin{equation}\label{eq:g_distinct}
		g(i, k) \neq g(j, m) \quad \text{for all } 1 \le j < i \le n-1, \ 1 \le k \le p_i-1, \text{ and } 1 \le m \le p_j-1.
	\end{equation}
\end{lemma}

\begin{proof}
	Suppose, for contradiction, that \eqref{eq:g_distinct} fails for some indices. Then $g(i, k) = g(j, m)$, which by Lemma~\ref{l5} is equivalent to
	\[
	\frac{P}{p_i} (p_i - k) = \frac{P}{p_j} (p_j - m) \iff k p_j = m p_i.
	\]
	Since $m p_i \equiv 0 \pmod{p_i}$, it follows that $k p_j \equiv 0 \pmod{p_i}$. By the coprimality of $p_i$ and $p_j$, we obtain $k \equiv 0 \pmod{p_i}$, which contradicts the condition $1 \le k \le p_i - 1$. \qedhere
\end{proof}

\begin{lemma} \label{l7}
The $T$-degrees of the vertices from $A_i$ are determined as follows:
\begin{enumerate}
\item $f(i, 1) = \dfrac{P}{2p_i} (S - p_i - n + 2) + p_{n-1}$ for all $1 \le i \le n-2$;
\item $f(i, k) = \dfrac{P}{2p_i} (S - p_i - n + 2k)$ for all $1 \le i \le n-2$ and $2 \le k \le p_i$;
\item $f(n-1, k) = \dfrac{P}{2p_{n-1}} (S - p_{n-1} - n + 2k) + 1$ for all $1 \le k \le p_{n-1}$.
\end{enumerate}
\end{lemma}
\begin{proof}
Fix $1 \le i \le n-1$ and $1 \le k \le p_i$. Let $h(i, k)$ denote the number of $T$-subgraphs of $D(T, \bar{p})$ containing the vertex $a_{i,k}$ but not containing $c$. Lemma~\ref{l3} implies that $h(i, k)$ equals the number of vertex subsets $\{v_1, v_2, \dots, v_n\}$ in $D(T, \bar{p})$ such that $v_j \in A_j$ for all $1 \le j \le n-1$, $v_i = a_{i,k}$, $v_n \in \bigcup_{l=1}^{n-1} B_l$, and $(v_j, v_n)^* \in A(D(T, \bar{p}))$ for all $1 \le j \le n-1$.
	
To determine $h(i, k)$, we compute the number of such subsets for all possible choices of $v_n \in B_l$ over all $1 \le l \le n-1$.
	
If $l = i$, the vertex $v_n \in B_i$ can be chosen in exactly $k-1$ ways, since only the vertices $b_{i,1}, \dots, b_{i,k-1}$ are adjacent to $a_{i,k}$ in $D(T, \bar{p})$. Moreover, these vertices are adjacent to all vertices in $A_j$ for all $j \neq i$, yielding $p_j$ choices for each $j$. Thus, the number of such subgraphs with $v_n \in B_i$ is
\[
(k-1) \prod_{\substack{j=1 \\ j \neq i}}^{n-1} p_j = \frac{P}{p_i} (k-1).
\]
	
	If $l \neq i$, fix an arbitrary vertex $v_n = b_{l,m} \in B_l$ with $1 \le m \le p_l-1$. Then $v_l \in A_l$ can be chosen in $p_l - m$ ways, as $b_{l,m}$ is adjacent to exactly $p_l - m$ vertices in $A_l$, namely $a_{l,m+1}, a_{l,m+2}, \dots, a_{l,p_l}$. For any $j \neq i, l$, the vertex $v_j \in A_j$ can be chosen in $p_j$ ways since $b_{l,m}$ is adjacent to all vertices in $A_j$. Hence, for a fixed $b_{l,m} \in B_l$, the number of choices is
	\[
	(p_l - m) \prod_{\substack{j=1 \\ j \neq i, l}}^{n-1} p_j = \frac{P}{p_i p_l} (p_l - m).
	\]
	Summing this expression over all $1 \le m \le p_l-1$ yields the total number of the required subsets with $v_n \in B_l$ ($l \neq i$):
	\[
	\sum_{m=1}^{p_l-1} \frac{P}{p_i p_l} (p_l - m) = \frac{P}{p_i p_l} \frac{(p_l-1)p_l}{2} = \frac{P}{2p_i} (p_l - 1).
	\]
Combining these cases, we evaluate $h(i, k)$ as follows:
\begin{align*}
	h(i, k) &= \sum_{\substack{l=1 \\ l \neq i}}^{n-1} \frac{P}{2p_i} (p_l - 1) + \frac{P}{p_i} (k-1) \\
	&= \frac{P}{2p_i} \sum_{l=1}^{n-1} (p_l - 1) - \frac{P}{2p_i} (p_i - 1) + \frac{P}{p_i} (k-1)
	= \frac{P}{2p_i} (S - p_i - n + 2k).
\end{align*}
			
The final equalities follow by combining $h(i, k)$ with the properties of the vertex $c$ in the graph $D(T, \bar{p})$ established in the proof of Lemma~\ref{l4}. 

For claim~1, fix $1 \le i \le n-2$. Since $a_{i,1}$ belongs to all $p_{n-1}$ $T$-subgraphs containing $c$, we obtain
\[
f(i, 1) = h(i, 1) + p_{n-1} = \frac{P}{2p_i} (S - p_i - n + 2) + p_{n-1}.
\]

For claim~2, fix $1 \le i \le n-2$ and $2 \le k \le p_i$. Because $a_{i,k}$ does not belong to any $T$-subgraph containing $c$, we have
\[
f(i, k) = h(i, k) = \frac{P}{2p_i} (S - p_i - n + 2k).
\]

For claim~3, fix $1 \le k \le p_{n-1}$. Since $a_{n-1,k}$ belongs to exactly one $T$-subgraph containing $c$, we have
\[
f(n-1, k) = h(n-1, k) + 1 = \frac{P}{2p_{n-1}} (S - p_{n-1} - n + 2k) + 1. \qedhere
\]
\end{proof}
		
\begin{corollary} \label{cor.3}
	$f(i, 1) < f(i, 2) < \dots < f(i, p_i)$ for all $1 \le i \le n-1$.
\end{corollary}
\begin{proof}
	The strict inequality $f(i, k+1) > f(i, k)$ for $1 \le k \le p_i - 1$ follows directly from the formulas of Lemma~\ref{l7} for all $k \ge 2$ (and for all $k \ge 1$ if $i = n-1$). In the remaining case where $1 \le i \le n-2$ and $k = 1$, the relation $f(i, 2) > f(i, 1)$ takes the form
	\[
	\frac{P}{2p_i} (S - p_i - n + 4) > \frac{P}{2p_i} (S - p_i - n + 2) + p_{n-1} \iff P > p_i p_{n-1},
	\]
	which holds by virtue of Lemma~\ref{l2}. 
\end{proof}

\begin{lemma} \label{l8}
	The $T$-degrees of any vertices from distinct sets $A_i$ and $A_j$ satisfy the inequality
	\begin{equation}\label{eq:f_distinct}
		f(i, k) < f(j, m) \quad \text{for all } 1 \le j < i \le n-1, \ 1 \le k \le p_i, \text{ and } 1 \le m \le p_j.
	\end{equation}
\end{lemma}
\begin{proof}
	Fix $1 \le j < i \le n-1$, $1 \le k \le p_i$, and $1 \le m \le p_j$. By Corollary~\ref{cor.3}, we have
	\[
	f(i, k) \le f(i, p_i) \quad \text{and} \quad f(j, 1) \le f(j, m).
	\]
	Therefore, to establish \eqref{eq:f_distinct}, it suffices to show that
	\begin{equation}\label{eq:f_blocks_max_min}
		f(i, p_i) < f(j, 1) \quad \text{for all } 1 \le j < i \le n-1.
	\end{equation}
	From Lemma~\ref{l7}, the following relations hold for all $1 \le j < i \le n-1$:
	\begin{align*}
		f(i, p_i) &\le \frac{P}{2p_i} (S - n + p_i) + 1, \\
		f(j, 1) &= \frac{P}{2p_j} (S - p_j - n + 2) + p_{n-1} > \frac{P}{2p_j} (S - p_j - n + 2) + 1.
	\end{align*}
	Thus, inequality \eqref{eq:f_blocks_max_min} holds provided that
	\begin{equation}\label{eq:f_sufficient_algebraic}
		\frac{P}{2p_j} (S - p_j - n + 2) > \frac{P}{2p_i} (S - n + p_i) \iff p_i (S - p_j - n + 2) > p_j (S - n + p_i).
	\end{equation}
	To prove \eqref{eq:f_sufficient_algebraic}, we apply the main condition \eqref{eq:main}, which yields $p_j > n$ and $p_i > 3p_j$ whenever $1 \le j < i \le n-1$. Combining these with the bound $S - n - 3p_j > 0$ from Lemma~\ref{l2}, we obtain:
	\begin{align*}
		p_i (S - p_j - n + 2) &> 3p_j (S - p_j - n + 2) > p_j (3S - 3p_j - 3n) \\
		&= p_j \big( (S - n + p_i) + (S - n - p_i) + (S - n - 3p_j) \big) \\
		&> p_j \big( (S - n + p_i) + (S - p_j - p_i) \big) > p_j (S - n + p_i),
	\end{align*}
	which completes the evaluation of \eqref{eq:f_sufficient_algebraic}. 
\end{proof}

The next property follows immediately from Lemma~\ref{l8} and Corollary~\ref{cor.3}.
\begin{corollary} \label{cor.4}
	The minimum $T$-degree among the vertices of the set $\bigcup_{i=1}^{n-1} A_i$ is $f(n-1, 1)$.
\end{corollary}

\begin{lemma} \label{l9}
	The $T$-degree of the vertex $c$ is strictly less than the $T$-degree of any vertex from the set $A_i$:
	\begin{equation}\label{eq:c_less_f}
		f(c) < f(i, k) \quad \text{for all } 1 \le i \le n-1 \text{ and } 1 \le k \le p_i.
	\end{equation}
\end{lemma}
\begin{proof}
	By Corollary~\ref{cor.4}, to establish \eqref{eq:c_less_f}, it suffices to show that $f(n-1, 1) > f(c)$. Applying Lemmas~\ref{l2}, \ref{l4}, \ref{l7}, and the bound $p_{n-2} > n$ from \eqref{eq:main}, we obtain:
	\begin{align*}
		f(n-1, 1) &= \frac{P}{2p_{n-1}} (S - p_{n-1} - n + 2) + 1 \\[1mm]
		&> \frac{P}{2p_{n-1}} \big( (p_{n-2} + p_{n-1}) - p_{n-1} - n + 2 \big)
		> \frac{P}{p_{n-1}} > p_{n-1} = f(c). \qedhere
		\end{align*}
		\end{proof}

\begin{lemma} \label{l10}
	The $T$-degree of the vertex $c$ is strictly less than the $T$-degree of any vertex from the set $B_i$:
	\begin{equation}\label{eq:c_less_g}
		f(c) < g(i, m) \quad \text{for all } 1 \le i \le n-1 \text{ and } 1 \le m \le p_i-1.
	\end{equation}
\end{lemma}
\begin{proof}
	By Lemmas~\ref{l4} and \ref{l5}, the inequality \eqref{eq:c_less_g} is equivalent to
	\[
	\frac{P}{p_i} (p_i - m) > p_{n-1} \quad \text{for all } 1 \le i \le n-1 \text{ and } 1 \le m \le p_i-1,
	\]
	the validity of which follows from Lemma~\ref{l2} and the estimate:
	\[
	\frac{P}{p_i} (p_i - m) \ge \frac{P}{p_i} > p_{n-1}. \qedhere
	\]
\end{proof}

\begin{lemma} \label{l11}
	The $T$-degrees of vertices from the sets $A_i$ and $B_i$ satisfy the following relations:
	\begin{enumerate}
		\item For all $1 \le i \le n-2$, $1 \le k \le p_i$, and $1 \le m \le p_i-1$:
		\begin{equation}\label{eq:f_greater_g}
			f(i, k) > g(i, m);
		\end{equation}
\item For all $1 \le k \le p_{n-1}$ and $1 \le m \le p_{n-1}-1$:
		\begin{equation}\label{eq:f_not_g}
			f(n-1, k) \neq g(n-1, m).
		\end{equation}
	\end{enumerate}
\end{lemma}
\begin{proof}
	For claim~1, fix $1 \le i \le n-2$, $1 \le k \le p_i$, and $1 \le m \le p_i-1$. By Lemmas~\ref{l5}, \ref{l7} and Corollaries~\ref{cor.2} and \ref{cor.3}, we have the following boundaries:
	\[
	f(i, k) \ge f(i, 1) = \frac{P}{2p_i} (S - p_i - n + 2) + p_{n-1} \quad \text{and} \quad g(i, m) \le g(i, 1) = \frac{P}{p_i} (p_i - 1).
	\]
	Since $p_{n-1} > 0$, relation \eqref{eq:f_greater_g} holds because the inequality
	\[
	\frac{P}{2p_i} (S - p_i - n + 2) > \frac{P}{p_i} (p_i - 1) \iff S > 3p_i + n - 4
	\]
	is valid strictly by virtue of Lemma~\ref{l2}.
	
	For claim~2, suppose, for contradiction, that \eqref{eq:f_not_g} fails, so that $f(n-1, k) = g(n-1, m)$. By Lemmas~\ref{l5} and \ref{l7}, this equality is equivalent to
	\begin{align*}
		\frac{P}{2p_{n-1}} (S - p_{n-1} - n + 2k) + 1 &= \frac{P}{p_{n-1}} (p_{n-1} - m) \\
		&\iff P(S - p_{n-1} - n + 2k) + 2p_{n-1} = 2P(p_{n-1} - m).
	\end{align*}
Since $p_1$ divides $P$, this relation implies that $p_1$ must also divide $2p_{n-1}$. However, for $n \ge 4$, this contradicts the fact that $p_1$ and $p_{n-1}$ are distinct prime numbers with $p_1 > 6 > 2$, thereby completing the proof.
\end{proof}

\begin{lemma} \label{l12}
	The $T$-degrees of any two vertices from the sets $A_i$ and $B_j$ with $i \neq j$, respectively, are unequal:
	\begin{equation}\label{eq:f_g_distinct}
		f(i, k) \neq g(j, m) \quad \text{for all } 1 \le i, j \le n-1, \ i \neq j, \ 1 \le k \le p_i, \text{ and } 1 \le m \le p_j-1.
	\end{equation}
\end{lemma}
\begin{proof}
	If $f(i, k) = g(j, m)$ for some $i \neq j$, we obtain a contradiction in two cases depending on $i$.
	
	For Case~1, let $i = n-1$. Then $1 \le j \le n-2$. By Lemmas~\ref{l5} and \ref{l7}, the equality becomes
	\begin{align}
		\frac{P}{2p_{n-1}}& (S - p_{n-1} - n + 2k) + 1 = \frac{P}{p_j} (p_j - m) \nonumber \\
		&\quad \iff p_j P(S - p_{n-1} - n + 2k) + 2p_{n-1} p_j = 2p_{n-1} P(p_j - m). \label{eq:f_g_distinct_collapse}
	\end{align}
	Since $P$ divides both $p_j P(S - p_{n-1} - n + 2k)$ and $2p_{n-1} P(p_j - m)$, \eqref{eq:f_g_distinct_collapse} implies that $P$ must also divide $2p_{n-1} p_j$, which forces $2p_{n-1} p_j \ge P$. However, Lemma~\ref{l2} yields
	\[
	P \ge p_{n-3} p_{n-2} p_{n-1} > 2 p_{n-2} p_{n-1} \ge 2 p_{n-1} p_j,
	\]
	where the strict inequality holds as $p_{n-3} \ge p_1 > 6$ by \eqref{eq:main}, and $p_j \le p_{n-2}$ since $j \le n-2$. This contradicts the condition $2p_{n-1} p_j \ge P$.
	
	For Case~2, let $1 \le i \le n-2$. In this case, by Lemmas~\ref{l2}, \ref{l5}, and the formulas of Lemma~\ref{l7}, we obtain the following bounds for $f(i, k)$ and $g(j, m)$:
	\begin{align*}
		f(i, k) &\ge \frac{P}{2p_i} (S - p_i - n + 2k) > \frac{P}{2p_i} \big( (3p_i + n) - p_i - n \big) = P, \\
		g(j, m) &= \frac{P}{p_j} (p_j - m) < P.
	\end{align*}
	These bounds yield $f(i, k) > g(j, m) $ for all $1 \le i \le n-2$ and $1 \le j \le n-1$ with $i \neq j$, which directly contradicts the assumption $f(i, k) = g(j, m)$.
\end{proof}

\begin{theorem} \label{t2}
	For any tournament $T$ of order $n \ge 4$ and any ordered tuple of prime numbers $\bar{p} = (p_1, p_2, \dots, p_{n-1})$ satisfying the main condition \eqref{eq:main}, the graph $D(T, \bar{p})$ is $T$-irregular.
\end{theorem}

\begin{proof}
	To establish the $T$-irregularity of $D(T, \bar{p})$, it suffices to verify that its vertices have pairwise distinct $T$-degrees. Based on the preceding results, the distribution of $T$-degrees is partitioned into non-overlapping classes:
	\begin{itemize}
		\item Within the sets $B_i$, the $T$-degrees strictly decrease inside each block by Corollary~\ref{cor.2} and are distinct for vertices from different blocks by Lemma~\ref{l6}.
		\item Within the sets $A_i$, the $T$-degrees strictly increase inside each block by Corollary~\ref{cor.3} and are distinct for vertices from different blocks by Lemma~\ref{l8}.
		\item The vertex $c$ has a strictly smaller $T$-degree than any vertex in $\bigcup_{i=1}^{n-1} A_i$ by Lemma~\ref{l9} and $\bigcup_{i=1}^{n-1} B_i$ by Lemma~\ref{l10}.
		\item The $T$-degrees of vertices from $A_i$ and $B_j$ are unequal, whether within the same block by Lemma~\ref{l11} or across distinct blocks by Lemma~\ref{l12}.
	\end{itemize}
	Consequently, all vertices in $V(D(T, \bar{p}))$ have pairwise distinct $T$-degrees, which implies that $D(T, \bar{p})$ is $T$-irregular.
\end{proof}

\smallskip
Theorem~\ref{t1} now follows from Proposition~\ref{p1}, Theorem~\ref{t2}, and Lemma~\ref{l1}.

\begin{corollary} \label{cor.5}
	For any tournament $T$ of order $n \ge 4$ and any graph $D(T, \bar{p})$, the underlying graph $U(D(T, \bar{p}))$ is $K_n$-irregular.
\end{corollary}
\begin{proof}
	By the proof of Lemma~\ref{l3}, a subgraph of $D(T, \bar{p})$ on $n$ vertices is a $T$-subgraph if and only if its vertices are pairwise adjacent. Since removing the arc orientations preserves vertex adjacency, the $K_n$-degree of a vertex $v$ in $U(D(T, \bar{p}))$ coincides with its $T$-degree in $D(T, \bar{p})$. Thus, the $T$-irregularity of $D(T, \bar{p})$ established in Theorem~\ref{t2} implies the $K_n$-irregularity of $U(D(T, \bar{p}))$.
\end{proof}


\begin{thebibliography}{99}
	\bibitem{1}{ G.~Chartrand, K.S.~Holbert,  O.R.~Oellermann, H.C.~Swart.}  $F$-degrees in graphs. {\it Ars Comb.}, 1987, vol.~24,  pp.~133--148.
	
	\bibitem{2}{A.~Ali, G.~Chartrand, P.~Zhang.} {\it Irregularity in graphs}, NY,  Springer, 2021, 109~p.\\
	\href{https://doi.org/10.1007/978-3-030-67993-4}{https://doi.org/10.1007/978-3-030-67993-4}
	
	\bibitem{3}{T.S.~Dovzhenok, A.V.~Filuta,  N.E.~Chuhai.} On some results of the study of \mbox{$F$-irregular} graphs in the class of biconnected graphs $F$. {\it Journal of the Belarus State University. Mathematics and Informatics}, 2024,  no.~2,  pp.~54--64 (in Russian).
	
	\bibitem{4}{T.S.~Dovzhenok.} Proof of the strong conjecture about $F$-irregular graphs in the class of graphs $\{F\}$ of diameter $2$. {\it Mat. Zametki}, 2026, vol.~120, no.~3, pp.~348--364. \href{https://doi.org/10.4213/mzm14573}{https://doi.org/10.4213/mzm14573}; {\it Math. Notes}, 2026, vol.~120, no.~3, pp.~439--452. \href{https://doi.org/10.1134/S0001434626603795}{https://doi.org/10.1134/S0001434626603795}
	
	\bibitem{5}{J.A.~Schreib.} {\it A construction of $F$-irregular graphs}, arXiv preprint arXiv:2609.07745, 2026.
	\href{https://doi.org/10.48550/arXiv.2609.07745}{https://doi.org/10.48550/arXiv.2609.07745}
	
	\bibitem{6}{T.~Dovzhenok, A.~Filuta.} {\it A Complete Proof of the Strong Conjecture about \mbox{$F$-Irregular} Graphs}, arXiv preprint arXiv:2609.19245, 2026.\\
	\href{https://doi.org/10.48550/arXiv.2609.19245}{https://doi.org/10.48550/arXiv.2609.19245}
	
	
	\bibitem{7}{T.~Dovzhenok.} {\it On the existence of $F$-irregular graphs for disconnected graphs $F$}, arXiv preprint arXiv:2610.18349, 2026. 
	\href{https://doi.org/10.48550/arXiv.2610.18349}{https://doi.org/10.48550/arXiv.2610.18349}
	
	\bibitem{8}{T.S.~Dovzhenok, I.D.~Lukashenko, Y.V.~Filiuta.} On $\overrightarrow{C_{n}}$-irregular oriented graphs. {\it Trudy Instituta Matematiki i Mekhaniki UrO RAN}, 2026. \\
	\href{https://doi.org/10.21538/0134-4889-2026-32-4-fon-01}{https://doi.org/10.21538/0134-4889-2026-32-4-fon-01}	 

			
	\end{thebibliography}
\end{document}